\documentclass[a4paper,12pt]{amsart}
\usepackage[utf8]{inputenc}
\usepackage{amsmath,amssymb,amsfonts,amsthm,graphicx}

\newtheorem{theorem}{Theorem}[section]

\numberwithin{equation}{section}

\theoremstyle{remark}

\newtheorem{definition}[theorem]{Definition}

\title{A family of periodic solutions of the three body problem. {\rm Light version}}

\author{Oscar Perdomo}
\thanks{Department of Mathematics, Central Connecticut State University, 120 Marcus White Hall, New Britain, CT 06052, USA \texttt{mihaib@ccsu.edu}}

\begin{document}

\begin{abstract}
In this paper we describe a 1-dimensional family of initial conditions $\Sigma$ that provides reduced periodic solution of the three body problem. This family $\Sigma$ contains a bifurcation point and extend the periodic solution described in \cite{P1}. This 1-dimensional family is the union of two embedded smooth curves. We will explain how the trajectories of the bodies in the solutions coming from one of the embedded curves have two symmetries while those coming from the other embedded curve only have one symmetry. The Round Taylor Method is a numerical method implemented by the author to keep track of the global error and the round-off error. A second version of this paper, same title with the ``light version'' part removed, will include analysis of the error of the solutions using the Round Taylor Method.
\end{abstract}

\maketitle

\begin{figure}[hbtp]
\begin{center}\includegraphics[width=.55\textwidth]{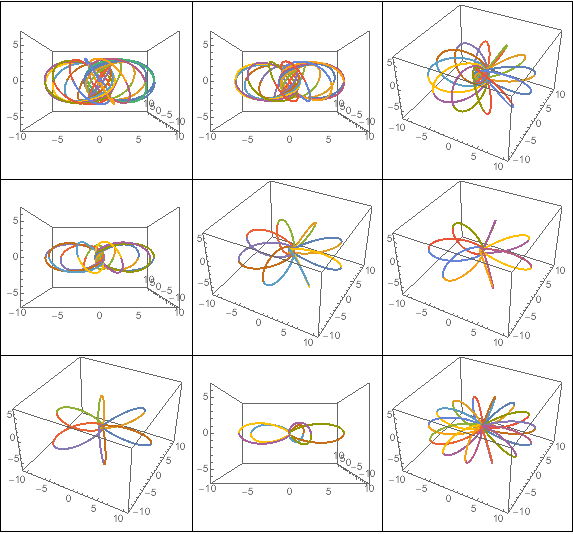}
\end{center}
\caption{Trajectory of one of the bodies for some of the periodic solutions described in this paper}\label{fig1}
\end{figure}

\section{Introduction}

In \cite{P1} the author considered the family of solutions of the three body problem described by  the parameters $a$ and $b$ suggested by the following figure.
\begin{figure}[hbtp]
\begin{center}\includegraphics[width=.5\textwidth]{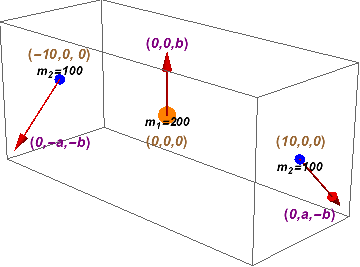}
\end{center}
\caption{A family of solutions of the three body problem described by the parameters $a$ and $b$. }\label{fig2}
\end{figure}
Three functions $F(t,a,b)$, $R(t,a,b)$ and $\Theta(t,a,b)$ describe the solution of the three body problem with initial conditions described in Figure \ref{fig2}, more precisely we have that

\begin{eqnarray*}\label{sol}
x(t)&=&(0,0,F(t,a,b))\\
y(t)&=&(R(t,a,b) \cos(\Theta(t,a,b)),R(t,a,b) \sin(\Theta(t,a,b)),- F(t,a,b))\\
z(t)&=&(-R(t,a,b) \cos(\Theta(t,a,b)),-R(t,a,b)\sin(\Theta(t,a,b)),-F(t,a,b))\, ,
\end{eqnarray*}

satisfy the three body problem equations, provided that the following ODE system holds true

\begin{eqnarray}\label{e1}
\ddot{F}=-\frac{400}{S^3}F,\, \quad \ddot{R}=\frac{100 a^2}{R^3}-\frac{25}{ R^2}-\frac{200 R}{S^3},
\, \quad R^2\dot{\Theta}=10 a,
\end{eqnarray}
with $\, R(0,a,b)=10, \, F(0,a,b)=0,\, \Theta(0,a,b)=0$, where  $S=\sqrt{R^2+4
F^2}$ 

where $\dot{G}$ denotes the partial derivative with respect to $t$ for any function $G(t,a,b)$. We are assuming that the gravitational constant is 1. We will be referring to the solution of the three body problem described by the differential equation  (\ref{sol})  as  $\phi(a,b)$.

where $\dot{G}$ denotes the partial derivative with respect to $t$ for any function $G(t,a,b)$. We are assuming that the gravitational constant is 1. We will be referring to the solution of the three body problem described by the differential equation  (\ref{sol})  as  $\phi(a,b)$.

It is not difficult to show that anytime we find a point $(T,a,b)$ such that $\dot{F}(T,a,b)=\dot{R}(T,a,b)=0$, then the solution  $\phi(a,b)$ is reduced periodic with period $4T$. We will call these solutions {\it odd/even} solutions because from the point of view of $t=0$ the function $F$ is an odd function with respect to $t$, but from the point of view of $t=T$, both functions  $F$ and $R$ are even. Also, it is not difficult to show that anytime we find a point $(T,a,b)$ such that $F(T,a,b)=\dot{R}(T,a,b)=0$, then the solution  $\phi(a,b)$ is reduced periodic with period $2T$. We will call these solutions {\it odd} solutions because from the point of view of $t=0$ the function $F$ is an odd function. We point out that every odd/even function is also odd due to the fact that if  $\dot{F}(T,a,b)=\dot{R}(T,a,b)=0$, then $F(2T,a,b)=\dot{R}(2T,a,b)=0$.
 
In \cite{P1} the author proved the existence of a small path $\alpha \subset R^3=\{(T,a,b): T, a ,b \in R\}$ of odd/even solutions. In this paper we extend this path to a path $S_1$ of odd/even solutions: the embedded curve $P_3BP_2$ in Figure \ref{fig3}. This path starts near a periodic solution where the body in the center remains motionless and the other two bodies move along an ellipse ($b$ near $0$). The path ends near a motion with a double coalition (in this case $a$ is near 0). We also provide a path $S_2\subset R^3$ of odd solutions: the embedded curve that contains the curve $P_1B$ in Figure \ref{fig3}. This path $S_2$ starts near a motion with a triple collision (in this case $a$ is near 0). This path $S_2$ intercepts the path $S_1$ in a single point that can be thought as a bifurcation point of the set of reduced periodic solutions. The author is not sure if this path will continue to be an unbounded curve. 

A second version of this paper will use the Round Taylor Method \cite{P1} to compute the error and show that the difference between the initial conditions and the values of the solution after a period $T$ is small.

\section{Graph of the path of initial conditions that provides reduced periodic solutions.}

\begin{figure}[hbtp]
\begin{center}\includegraphics[width=.5\textwidth]{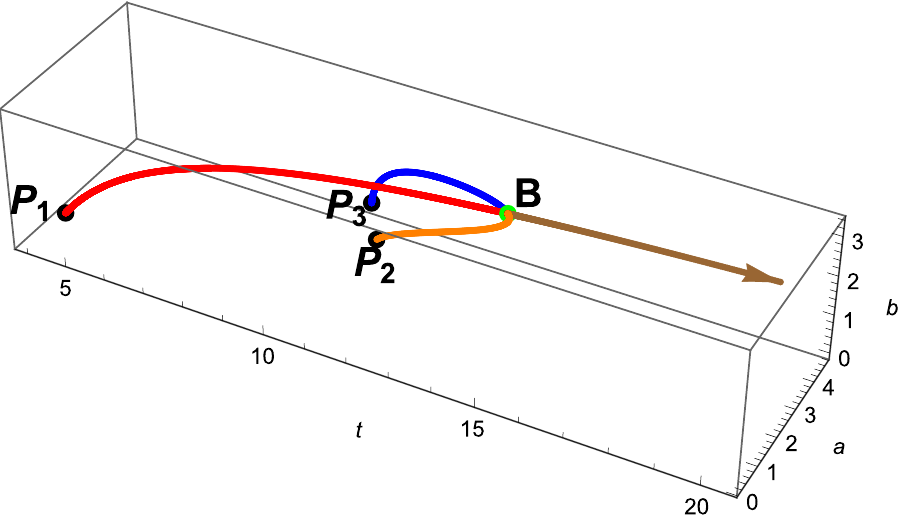}
\end{center}
\caption{$P_1$ has coordinates $(5.0063\dots, 0.109392\dots, 1.16473\dots)$, which means that when $a=0.109392\dots$, and $b=1.16473\dots$, then, the solution $\phi(a,b)$ has period $T=5.0063\dots$. The other coordinates are: $P_2=(12.7012\dots, 0.0437163\dots, 3.19541\dots)$, $P_3=(9.9472\dots, 4.73605\dots, 0.2)$ and,  the coordinate for bifurcation point are $B=(14.6072\dots, 2.08181\dots , 3.19493\dots)$.}\label{fig3}
\end{figure}

This section describes a path $\subset R^3=\{(T,a,b): T, a ,b \in R\}$ with the property that for any 
$(T,a,b)$ in $S$, the  solution $\phi(a,b)$ is reduced periodic with period $T$, this is, the relative position of the three bodies in the solution $\phi(a,b)$ repeats every $T$ units of time.  Figure \ref{fig3} is an image of the path $S$.

Recall that the distances between the three bodies depend exclusively only on the functions $R(t)$ and $F(t)$. More precisely, the distance between the two bodies that go around the $z$ axis is $2 R(t)$ and the distance between the body that moves on the $z$-axis and any of the other two bodies is given by $\sqrt{R(t)^2+F(t)^2}$.


\subsection{The solution given by the point $P_3=(9.9472\dots, 4.73605\dots, 0.2)$} 
\begin{figure}[hbtp]
\begin{center}\includegraphics[width=.3\textwidth]{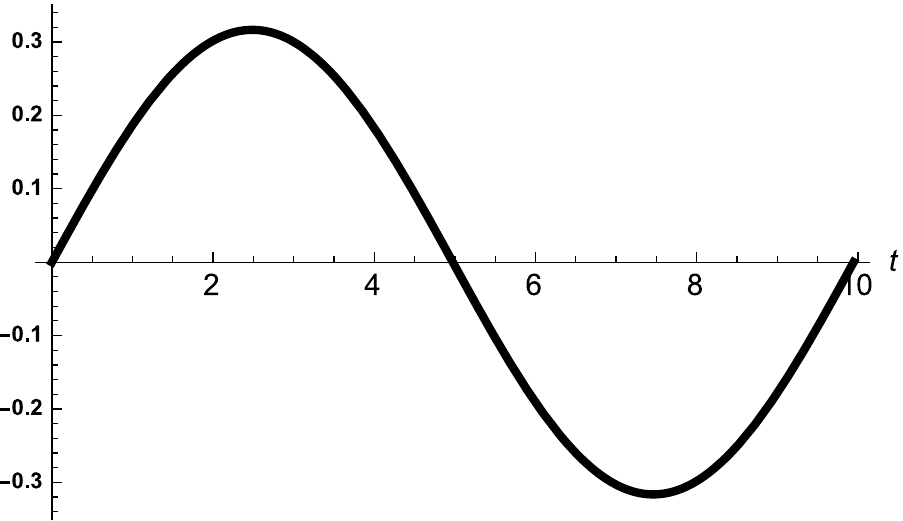}\hskip.2cm \includegraphics[width=.3\textwidth]{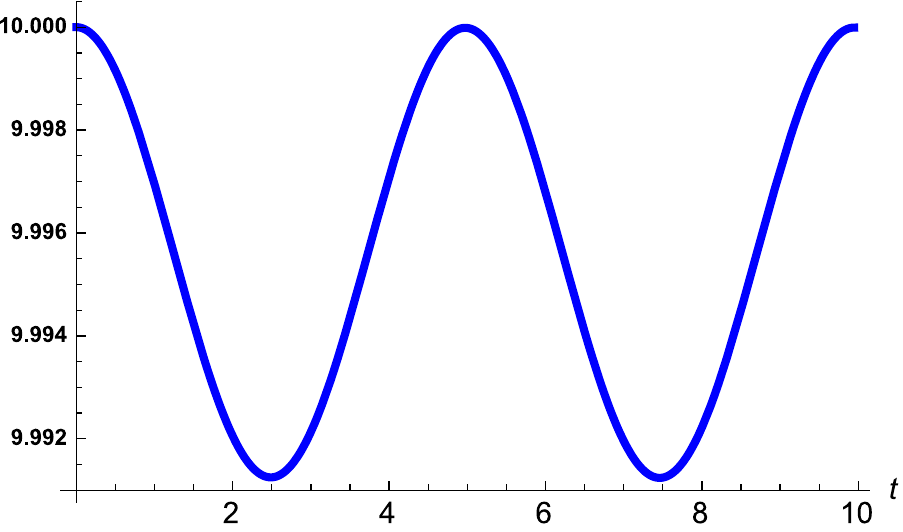}\hskip.2cm \includegraphics[width=.3\textwidth]{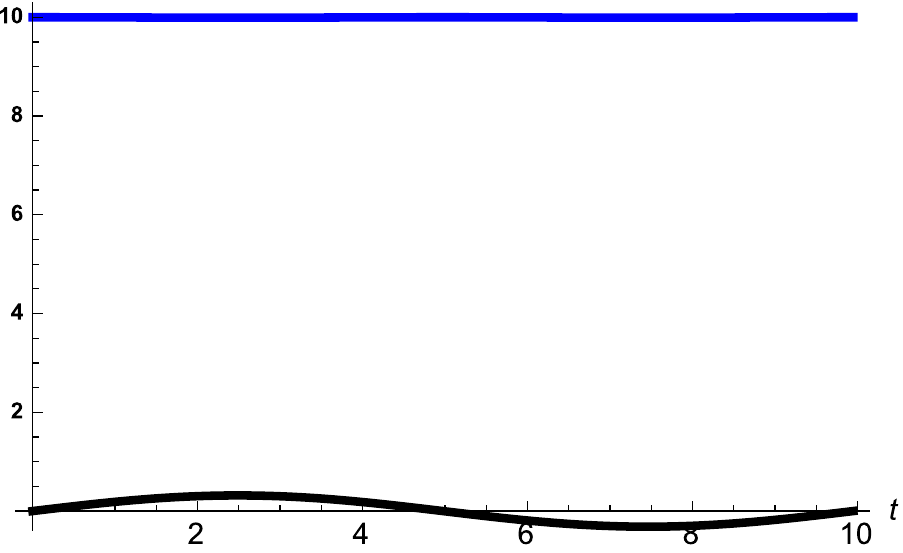}
\end{center}
\caption{For the point $P_3$, the image on the left shows the graph of the function $F$, the image on the center shows the graph of the function $R$ and the image on the right shows both functions.}\label{fig4}
\end{figure}

This solution is very close to the solution when the body in the center stays still and the other two bodies move on a perfect circle. Figure \ref{fig4} shows the graphs of the functions $F$ and $R$ associated with this solution. We can see how the body in the center given by $(0,0,F(t))$ moves very little and the other two bodies given by $(R(t) \cos(\Theta(t)),R(t) \sin(\Theta(t)),- F(t))$ and $(-R(t) \cos(\Theta(t)),-R(t)\sin(\Theta(t)),-F(t))$ stay near the circle of radius 10 on the $x$-$y$ plane with center at the origin. 

\subsection{The solution given by the point $P_2=(12.7012\dots, 0.0437163\dots, 3.19541\dots)$} 
\begin{figure}[hbtp]
\begin{center}\includegraphics[width=.5\textwidth]{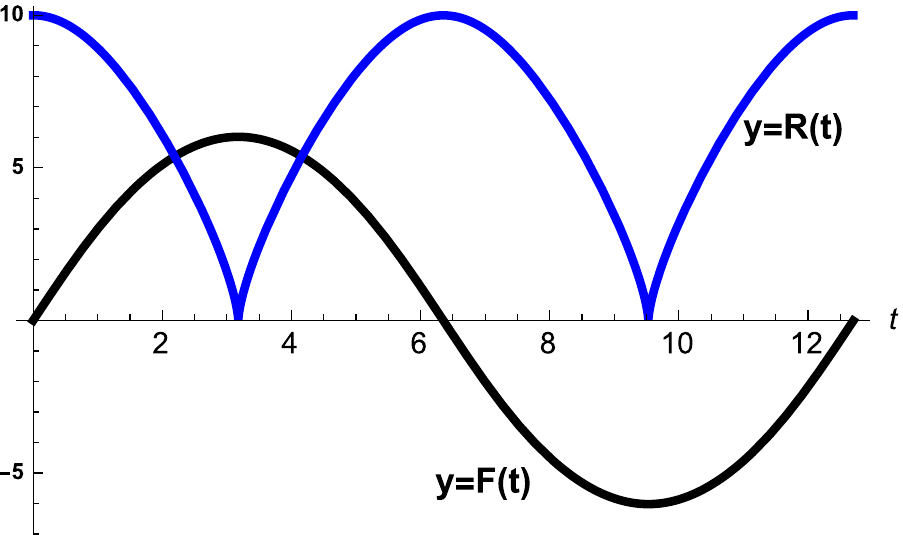}
\end{center}
\caption{The functions $F$ and $R$ associated with the reduced periodic solution given by the point $P_2$.}\label{fig5}
\end{figure}

For this motion the body in the center oscillates from $(0,0,-6.01695\dots)$ to $(0,0,6.01695\dots)$. This solution is very close to a solution with a double collision. Figure \ref{fig5} shows the graphs of the functions $F$ and $R$ associated with this solution. We can see how after a quarter of a period the two bodies that go around the $z$-axis are very close to each other ($R$ is very small, it is near $0.0547972$) at this instance the body in the center is at its highest point while the other two bodies are both very close to the point $(0,0,-6.01695\dots)$.

\subsection{The solution given by the point $P_1=(5.0063\dots, 0.109392\dots, 1.16473\dots)$} 
\begin{figure}[hbtp]
\begin{center}\includegraphics[width=.5\textwidth]{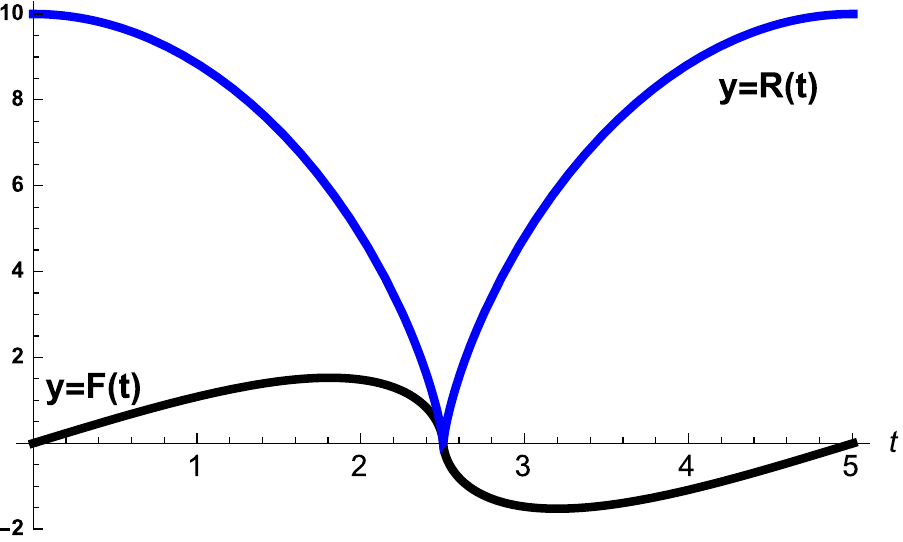}
\end{center}
\caption{The functions $F$ and $R$ associated with the reduced periodic solution given by the point $P_1$.}\label{fig6}
\end{figure}

For this motion the body in the center oscillates from $(0,0,-1.52497\dots)$ to $(0,0,1.52497\dots)$. This solution is very close to a solution with a triple collision. Figure \ref{fig6} shows the graphs of the functions $F$ and $R$ associated with this solution. We can see how after half of a period the three bodies are very close to the origin. $F$ is zero and $R$ is very small, it is near $0.0369262$.

\subsection{The solution given by the point $B=(14.6072\dots, 2.08181\dots , 3.19493\dots)$} \label{dp} 
\begin{figure}[hbtp]
\begin{center}\includegraphics[width=.5\textwidth]{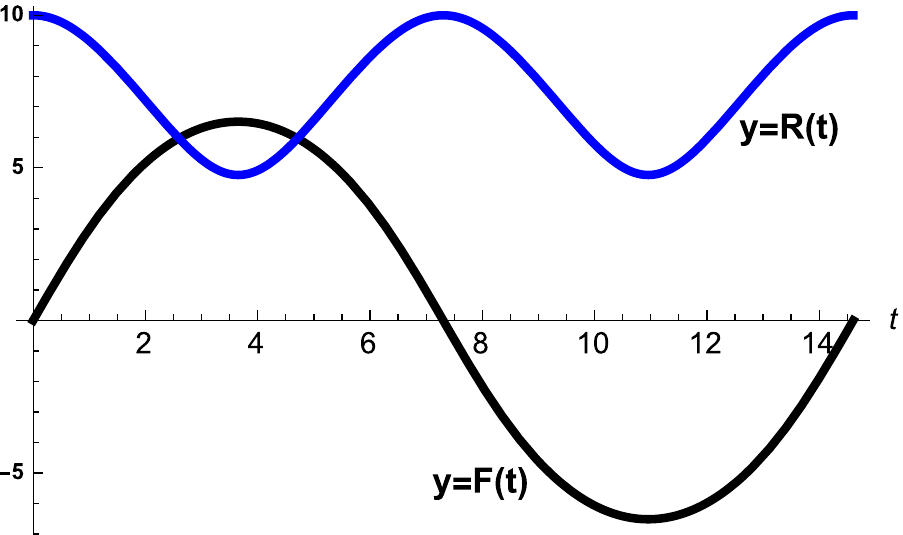}
\end{center}
\caption{The functions $F$ and $R$ associated with the reduced periodic solution given by the point $B$}\label{fig7}
\end{figure}

This is the bifurcation point of the path of solutions.  Figure \ref{fig7} shows the graphs of the functions $F$ and $R$ associated with this solution.

\begin{figure}[hbtp]
\begin{center}\includegraphics[width=.3\textwidth]{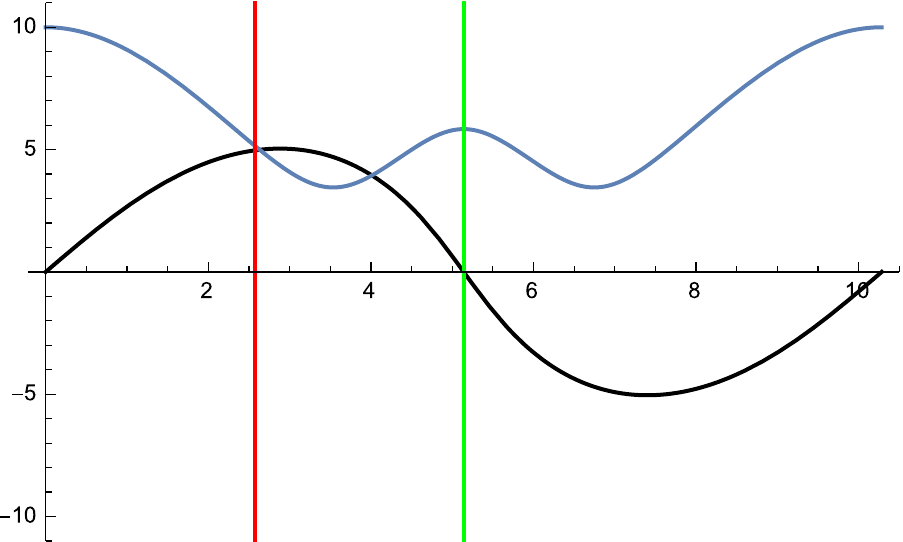}\hskip.6cm \includegraphics[width=.3\textwidth]{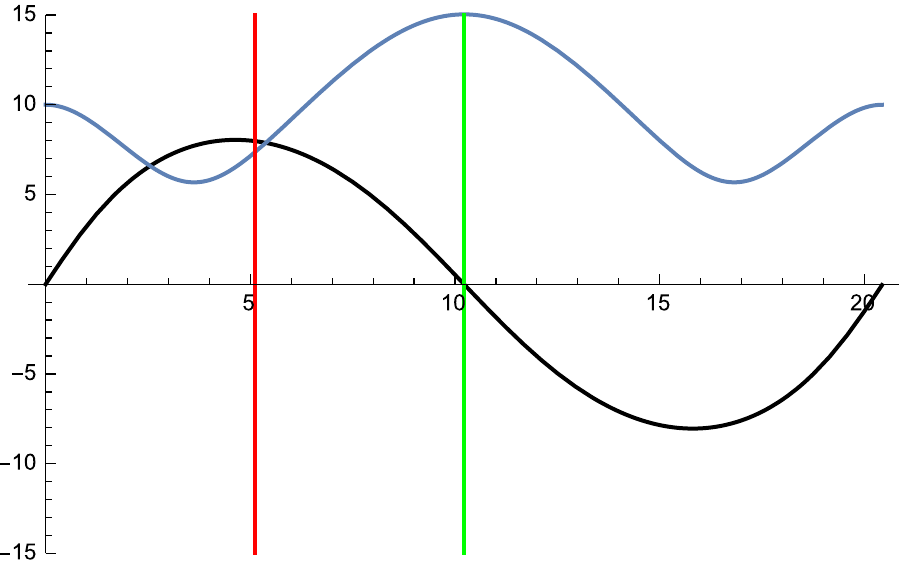}
\end{center}
\caption{The image on the left shows the functions $F$ and $R$ associated with a point on the curve $S_2$ between $P_1$ and $B$. The image on the right shows these two functions for a point on the curve $S_2$ located after the point $B$. We can see how the symmetry with respect to the vertical line $t=\frac{T}{2} $ that is present for solutions coming from points in $S_1$ is not present in these solutions.}\label{fig8}
\end{figure}

 Let us call $S_1$, the part of the path contained in $S$ given by the smooth curve that connects $P_2$, $B$ and $P_3$, see Figure \ref{fig3}. Likewise, let us call $S_2$, the part of the path contained in $S$ given by the smooth curve that connects $P_1$, $B$ and also contains the part of $S$ that has an arrow at the end. All the points in the path $S_1$ are {\it odd/even} solutions; all of them share the symmetry given by the functions shown in Figure \ref{fig7}. More precisely: With respect to the origin, the function $F$ is odd and the function $R$ is even and with respect to $t=\frac{T}{4}$, both functions $F$ and $R$ are even. The solutions associated with points in the curve $S_2$ different from the point $B$ do not share these two symmetries. With respect to $t=0$, the function $F$ is odd and the function $R$ is even. There is not  symmetry with respect to $t=\frac{T}{4}$. Figure \ref{fig8} shows the function $F$ and $R$ for the solution associated with two points on $S_2$.
 
\section{Procedure to obtain the points on the curve $S$ and the bifurcation point}

\subsection{Original ODE} The functions $F(t,a,b)$,  $R(t,a,b)$ and $\Theta(t,a,b)$ satisfy the following ODE

\begin{eqnarray*}
\dot{x}_1&=& x_3 \cr
\dot{x}_2& =& x_4\cr
\dot{x}_3& =&-\frac{400 x_1}{\left(4 x_1^2+x_2^2\right){}^{3/2}}\cr
\dot{x}_4& =& \frac{100 a^2}{x_2^3}-\frac{25}{x_2^2}-\frac{200 x_2}{\left(4 x_1^2+x_2^2\right){}^{3/2}}\cr
\dot{x}_5&=& \frac{10 a}{x_2^2}
\end{eqnarray*}

with initial conditions

$$x_1(0)=0\quad x_2(0)=10\quad x_3(0)=b\quad x_4(0)=0\quad x_5(0)=0  $$ 

This ODE will be referred to as the {\it original ODE}. We have that $x_1(t)$ provides $F(t,a,b)$, $x_2(t)$ provides $R(t,a,b)$, $x_3(t)$ provides $\dot{F}(t,a,b)$, $x_4(t)$ provides $\dot{R}(t,a,b)$ and  $x_5(t)$ provides $\Theta(t,a,b)$

\subsection{Extended ODE} \label{odes} The partial derivative of the functions $F(t,a,b)$,  $R(t,a,b)$ and $\Theta(t,a,b)$ satisfy the following ODE

\begin{eqnarray*}
\dot{x}_1&=&x_3\cr 
\dot{x}_2&=&x_4\cr 
\dot{x}_3&=&-\frac{400 x_1}{\left(4 x_1{}^2+x_2{}^2\right){}^{3/2}}\cr 
\dot{x}_4&=&-\frac{200 x_2}{\left(4 x_1{}^2+x_2{}^2\right){}^{3/2}}-\frac{25}{x_2{}^2}+ \frac{100 a^2}{x_2{}^3}\cr 
\dot{x}_5&=&\frac{10 a}{x_2{}^2}\cr 
\dot{x}_6&=&x_8\cr 
\dot{x}_7&=&x_9\cr 
\dot{x}_8&=&\frac{3200 x_6 x_1{}^2}{\left(4 x_1{}^2+x_2{}^2\right){}^{5/2}}+\frac{1200 x_2 x_7 x_1}{\left(4 x_1{}^2+x_2{}^2\right){}^{5/2}}-\frac{400 x_2{}^2 x_6}{\left(4 x_1{}^2+x_2{}^2\right){}^{5/2}}\cr 
\dot{x}_9&=&-\frac{300 a^2 x_7}{x_2{}^4}+\frac{200 a}{x_2{}^3}+\frac{2400 x_1 x_2 x_6}{\left(4 x_1{}^2+x_2{}^2\right){}^{5/2}}+\frac{50 x_7}{x_2{}^3}+\frac{600 x_2{}^2 x_7}{\left(4 x_1{}^2+x_2{}^2\right){}^{5/2}}-\frac{200 x_7}{\left(4 x_1{}^2+x_2{}^2\right){}^{3/2}}\cr 
\dot{x}_{10}&=&\frac{10}{x_2{}^2}-\frac{20 a x_7}{x_2{}^3}\cr 
\dot{x}_{11}&=&x_{13}\cr
\dot{x}_{12}&=&x_{14}\cr 
\dot{x}_{13}&=&\frac{3200 x_{11} x_1{}^2}{\left(4 x_1{}^2+x_2{}^2\right){}^{5/2}}+\frac{1200 x_2 x_{12} x_1}{\left(4 x_1{}^2+x_2{}^2\right){}^{5/2}}-\frac{400 x_2{}^2 x_{11}}{\left(4 x_1{}^2+x_2{}^2\right){}^{5/2}}\cr 
\dot{x}_{14}&=&-\frac{300 a^2 x_{12}}{x_2{}^4}+\frac{2400 x_1 x_2 x_{11}}{\left(4 x_1{}^2+x_2{}^2\right){}^{5/2}}+\frac{50 x_{12}}{x_2{}^3}+\frac{600 x_2{}^2 x_{12}}{\left(4 x_1{}^2+x_2{}^2\right){}^{5/2}}-\frac{200 x_{12}}{\left(4 x_1{}^2+x_2{}^2\right){}^{3/2}}\cr 
\dot{x}_{15}&=&-\frac{20 a x_{12}}{x_2{}^3}
\end{eqnarray*}

with initial conditions
\begin{eqnarray*}
x_1(0)&=&0\cr x_2(0)&=&10\cr x_3(0)&=&b\cr x_4(0)&=&0\cr x_5(0)&=&0\cr x_6(0)&=&0\cr x_7(0)&=&0\cr x_8(0)&=&0\cr x_9(0)&=&0\cr x_{10}(0)&=&0\cr x_{11}(0)&=&0\cr x_{12}(0)&=&0\cr x_{13}(0)&=&1\cr x_{14}(0)&=&0\cr x_{15}(0)&=&0
\end{eqnarray*}

This ODE will be referred to as the {\it extendded ODE}. We have that $x_1(t)$ provides $F(t,a,b)$, $x_2(t)$ provides $R(t,a,b)$, $x_3(t)$ provides $\dot{F}(t,a,b)$, $x_4(t)$ provides $\dot{R}(t,a,b)$, $x_5(t)$ provides $\Theta(t,a,b)$, $x_6(t)$ provides $\frac{\partial F}{\partial a}(t,a,b)$, $x_7(t)$ provides $\frac{\partial R}{\partial a}(t,a,b)$, $x_8(t)$ provides $\frac{\partial^2 F}{\partial a\partial t}(t,a,b)$, $x_9(t)$ provides $\frac{\partial^2 R}{\partial a\partial t}(t,a,b)$,  $x_{10}(t)$ provides $\frac{\partial \Theta}{\partial a}(t,a,b)$, $x_{11}(t)$ provides $\frac{\partial F}{\partial b}(t,a,b)$, $x_{12}(t)$ provides $\frac{\partial R}{\partial b}(t,a,b)$, $x_{13}(t)$ provides $\frac{\partial^2 F}{\partial b\partial t}(t,a,b)$, $x_{14}(t)$ provides $\frac{\partial^2 R}{\partial b\partial t}(t,a,b)$ and  $x_{15}(t)$ provides $\frac{\partial \Theta}{\partial b}(t,a,b)$.

These equations easily follows as an application of the chain rule. To exemplify this procedure let us compute $\dot{x}_8(t)$,

\begin{eqnarray*}
\dot{x}_8 &=& \ddot{x}_1{}_a = \left( -400 x_1\, \left(4 x_1{}^2+x_2{}^2\right){}^{-3/2}\ \right)_a\cr
&= & -400 x_6\, \left(4 x_1{}^2+x_2{}^2\right){}^{-3/2} + 1200 x_1 \left(4 x_1{}^2+x_2{}^2\right){}^{-5/2} \, \left( 4 x_1 x_6+ x_2 x_7  \right)\cr
&= &\frac{3200 x_6 x_1{}^2}{\left(4 x_1{}^2+x_2{}^2\right){}^{5/2}}+\frac{1200 x_2 x_7 x_1}{\left(4 x_1{}^2+x_2{}^2\right){}^{5/2}}-\frac{400 x_2{}^2 x_6}{\left(4 x_1{}^2+x_2{}^2\right){}^{5/2}}
\end{eqnarray*}

\subsection{Getting the points on $S_1$}

Ideally we want each point $(t,a,b)\in S_1$ to satisfy the system of equations $\{\dot{F}(\frac{t}{4},a,b)=\dot{R}(\frac{t}{4},a,b)=0\}$.  
From the paper \cite{P1} we know that there exist a point $\bar{P}_0=(\bar{t},\bar{a},\bar{b})$ near 

$$P_0=(t_0,a_0,b_0)=\left(\frac{13366894627923}{5000000000000},\frac{43170475352787}{10000000000000},\frac{1490359743}{1000000000} \right)$$

such that $\dot{F}(\bar{P}_0)=\dot{R}(\bar{P}_0)=0$. If we assume for a moment that we completely know the functions $F$ and $R$ in the whole ${ R}^3$, then the way to find the curve $S_1$ would be simply, first, find the vector field $X=v_1\times v_2$, where $v_1=\nabla \dot{F}=(\ddot{F},\dot{F}_a,\dot{F}_b)$ and $v_2=\nabla \dot{R}=(\ddot{R},\dot{R}_a,\dot{R}_b)$ and second, find the integral curve of the vector field  $X$  that goes through $\bar{P}_0$. Recall that the first entry of points in $S_1$ is 4 times the first entry of points in this integral curve. Notice that using the notation from section \ref{odes} we have that

$$\nabla \dot{F}=\left(-\frac{400 x_1(t)}{\left(4 x_1(t){}^2+x_2(t){}^2\right){}^{3/2}},x_{8},x_{13}\right) $$

and

$$\nabla \dot{R}=\left(\frac{25 \left(-\frac{8 x_2(t){}^4}{\left(4 x_1(t){}^2+x_2(t){}^2\right){}^{3/2}}-x_2(t)+4 a^2\right)}{x_2(t){}^3}, x_{9},x_{14}\right) $$

Even though we do not know the vector field $\nabla\dot{F}$ and $\nabla\dot{R}$ everywhere, we can still use the idea of integrating the vector field $X$ using the Euler method. The next subsection explains the algorithm that we are using in this paper to get the points in $S_1$.

\subsubsection{Continuation algorithm}\label{alg}
For the Algorithm we select three small numbers $\epsilon_1$ and $\epsilon_2$ and $\epsilon_3$ that we use as tolerance for the error.  We want $0<\epsilon_1<\epsilon_2$. Recall that we want to find solutions of the system $\{\dot{F}=\dot{R}=0\}$. We will be collecting the solution of this system of equations in a set called $TS$. To start with, we make, $TS=\{P_0\}$, the set which only element is the numerical solution that we know. We will use $P_0$ to start the algorithm.

\begin{enumerate}
\item
Consider a point $Q_0=(\tilde{t},\tilde{a},\tilde{b})$ such that $|\dot{F}(Q_0)|<\epsilon_1$ and $|\dot{R}(Q_0)|<\epsilon_1$.

\item
Make $Y=Q_0$ and select a positive real number $h$ and a positive integer  $k$.

\item For $i=1$ to $i=k$ Do
\begin{itemize}
\item
Find the values of $X$ by numerically solving the extended ode (see section \ref{odes}) using the  values of $t$, $a$ and $b$ given by the entries of  $Y$. This is, we integrate for $t$ units of time where $t$ is the first entry of $Y$ and we use the values of $a$ and $b$ given by the second and third entry of $Y$ to provide the initial conditions of the ode.
\item
Make  $Y=Y+\frac{h}{|X|}\, X$ and $q_i=Y$.
\end{itemize}
\item
If for all $i=1\dots k$, $|\dot{F}(q_i)|<\epsilon_2$ and $|\dot{R}(q_i)|<\epsilon_2$ then go to the next numeral, otherwise go back to numeral (2) and select a smaller $h$ and/or a smaller $k$.

\item Find $Q_1$ such that $|\dot{F}(Q_1)|<\epsilon_1$ and $|\dot{R}(Q_1)|<\epsilon_1$ and such that it is near $q_k$ in the sense that $|q_k-Q_1|<\epsilon_3$

\item Make $TS=TS\cup \{Q_1\}\cup\cup_{k=1}^k\{q_i\}$

\item Go back to numeral (1) to start the process over by making $Q_0=Q_1$.

\end{enumerate} 
 
\begin{definition}
We will call \underline{pillar points} all the points $Q_i$ in the previous algorithm.
\end{definition} 

In our case the algorithm stopped or became difficult to carry near the point $P_3$ and $P_2$ because the ode was approaching a singularity near $P_2$ (recall that $P_2$ is close to a collision) and on the other hand the vector field $X$ is close to the zero vector near $P_3$.

As an example of the method, if we start with $Q_0=P_0$ and we use $\epsilon_1=10^{-6}$, $\epsilon_2=\epsilon_3=0.00005$,  $h=0.001$  and $k=200$, then, 

$q_{200}= (2.7219062659312807, 4.212655007080421, 1.6538674269975053)$ 

and we can use  $Q_1=(2.72191575576588, 4.212633490447383, 1.6538497779324066)$.

 \subsection{Getting the points on $S_1$}


Ideally we want each point $(t,a,b)\in S_2$ to satisfy the system of equations $\{F(\frac{t}{2},a,b)=\dot{R}(\frac{t}{2},a,b)=0\}$.  
If we assume for a moment that we completely know the functions $F$ and $R$ in the whole ${ R}^3$, then the way to find the curve $S_2$ would be simply, first, find the vector field $Z=w_1\times w_2$, where $w_1=\nabla F=(\dot{F},F_a,F_b)$ and $w_2=\nabla \dot{R}=(\ddot{R},\dot{R}_a,\dot{R}_b)$ and second, find the integral curve of the vector field  $Z$.  Recall that the first entry of points in $S_2$ is twice the first entry of points in this integral curve. Notice that using the notation from section \ref{odes} we have that $\nabla F=\left(x_3,x_6,x_{11}\right) $. The points in $S_2$ are found using the algorithm described in section \ref{alg}.

\subsection{Getting the bifurcation point.}

Recall that the solutions associated with points on $S_1$ are called odd/even solutions and those associated with points in $S_2$ are called even solutions. Having in mind the symmetries of the solutions, it is not a surprise that for points on the integral curve $\gamma$ of the vector field $X$ that passes to the point $\bar{P}_0$, the vector field $Z$ is parallel to the vector field $X$. The bifurcation point occurs because along $\gamma$, there is a point where the vector field $Z$ vanishes. In order to numerically compute the bifurcation point, we take the minimum of the set $ \{|Z(\frac{t}{2},a,b)|:(t,a,b)\in S_1\} $. This minimum occurs at the point $B=(14.607249047056753, 2.081806260749908, 3.194934273913219) \in S_1$.

\section{Periodic Solutions}

Recall that  a reduced periodic solution $\phi(a,b)$ with period $T$ is periodic if $\frac{\Theta(T,a,b)}{\pi}$ is a rational number. It is not difficult to see the values of the function $\Theta$ on the points on the $S$. In this section we will be selecting some values for $\Theta$ along the path $S$ that provides periodic solutions.  

\subsection{Periodic solution displayed on Figure \ref{fig1}}

Every entry of the matrix  in  Figure \ref{fig1} (located at the begging of the paper) shows the trajectory  of one of the three bodies that moves periodically solving the three-body problem. The following table provides the initial conditions of these 9 periodic solutions. The first entry $T$ represent the period of the reduced periodic solution. This is, the positions and velocities after  $T$ units of time agree with the initial positions and velocities up to a rotation of $\theta$ radians. This rotation $\theta$ is given in the second column of the table.  Every color in these images show the trajectory after $T$ units of time. These points are part of the path $S_2$. These periodic solutions only one symmetry. The periodic solution given by the last row in this table is very close to a triple collision.

\begin{center}
    \begin{tabular}{ | l | l | l |}
    \hline
    Coordinates:\hskip .4cm $(T,a,b)$ & $\Theta(T,a,b)$ & Image \\
    & &  location \\ \hline
    (7.464725167070125, 1.332130235206886, 2.39131806234605)&$\frac{14 \pi }{9}$ & (1,1)  \\ \hline
    (6.692611549615348, 1.1189865713077587, 2.178542324667063)& $\frac{11 \pi }{7} $& (1,2)  \\ \hline
   (6.156645197408822, 0.9203809141218081, 1.9818897189116862)&$ \frac{8 \pi }{5}$ & (1,3) \\ \hline 
     (5.758992220584509, 0.72863210825991, 1.791870916952545)& $\frac{13 \pi }{8}$ & (2,1) \\
    \hline  (5.454531817007846, 0.5397530375918577, 1.6029431700479464) &$ \frac{5 \pi }{3}$ & (2,2)\\
    \hline (5.32949878269853, 0.4457220498597047, 1.5077112808284203) & $\frac{12 \pi }{7}$& (2,3)\\
    \hline(5.220379172126002, 0.35201283725645105, 1.411839280497764)&$ \frac{7 \pi }{4}$&(3,1) \\
    \hline (5.1264630045948305, 0.25902392732426605, 1.3159094056731593) & $\frac{16 \pi }{9}$& (3,2) \\
    \hline 
    (5.096931182326409, 0.22679159236240487, 1.2822928944673215) & $\frac{9 \pi }{5}$& (3,3) \\
    \hline 
        \end{tabular}
\end{center}

\begin{figure}[hbtp]
\begin{center}\includegraphics[width=1.1\textwidth]{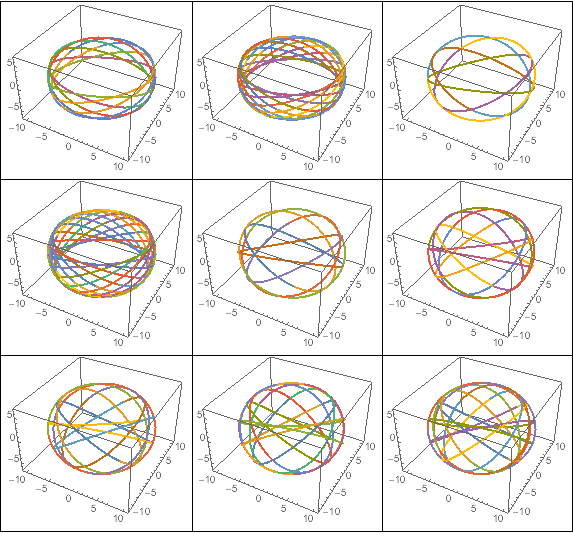}
\end{center}
\caption{Trajectory of one of the bodies for some of the periodic solutions described in this paper.}\label{f9}
\end{figure}

\subsection{Periodic solution displayed on Figure \ref{f9}}

Every entry of the matrix  in  Figure \ref{f9} (located at the end of the paper) shows the trajectory  of one of the three bodies that moves periodically solving the three-body problem. The following table provides the initial conditions of these 9 periodic solutions. The first entry $T$ represent the period of the reduced periodic solution. This is, the positions and velocities after  $T$ units of time agree with the initial positions and velocities up to a rotation of $\theta$ radians. This rotation $\theta$ is given in the second column of the table.  Every color in these images show the trajectory after $T$ units of time. These points are part of the path $S_1$. These periodic solutions are odd/even solutions and therefore they have two symmetries.

\begin{center}
    \begin{tabular}{ | l | l | l |}
    \hline
    Coordinates:\hskip .4cm $(T,a,b)$  & $\Theta(T,a,b)$ & Image \\
    & &  location \\ \hline
    (10.694782129146047, 4.3162773916465715, 1.4916623030663265)&$ \frac{17 \pi }{8}$ & (1,1)  \\ \hline
    (10.900907564917922, 4.205530359018152, 1.6642665366638942)& $\frac{15 \pi }{7} $& (1,2)  \\ \hline
   (11.251546754899636, 4.020933846016405, 1.9098147447184282)&$ \frac{13 \pi }{6}$ & (1,3) \\ \hline 
     (11.546587626864484, 3.8684643701482133, 2.0840362602898437)& $\frac{11 \pi }{5}$ & (2,1) \\
    \hline  (12.005866456188725, 3.6348152391561275, 2.314439960284758) &$ \frac{9 \pi }{4}$ & (2,2)\\
    \hline (12.482298481771874, 3.3941653076488447, 2.5161535540867117) & $\frac{16 \pi }{7}$& (2,3)\\
    \hline(12.805950195048542, 3.229502929775678, 2.6373059506494907)&$ \frac{7 \pi }{3}$&(3,1) \\
    \hline (13.037980888481863, 3.10968118474065, 2.717938296845948) & $\frac{12 \pi }{5}$& (3,2) \\
    \hline 
    (13.211458502729283`, 3.01851140163578`, 2.775360523910528`) & $\frac{17 \pi }{7}$& (3,3) \\
    \hline 
        \end{tabular}
\end{center}

\begin{figure}[hbtp]
\begin{center}\includegraphics[width=1.1\textwidth]{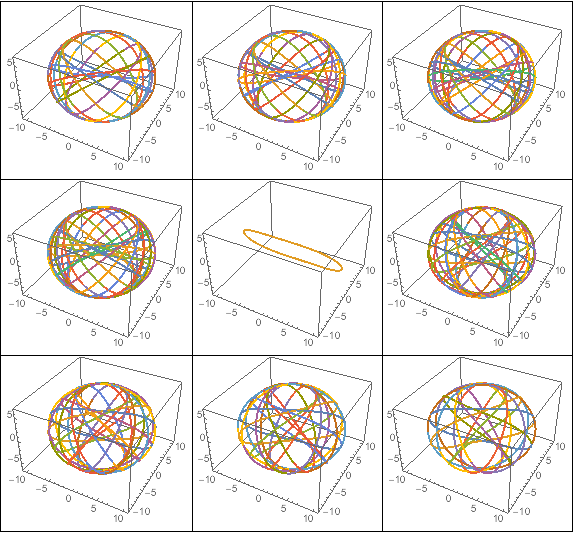}
\end{center}
\caption{Trajectory of one of the bodies for some of the periodic solutions described in this paper.}\label{fs9}
\end{figure}

\subsection{Periodic solution displayed on Figure \ref{fs9}}

Every entry of the matrix  in  Figure \ref{fs9} (located at the end of the paper) shows the trajectory  of one of the three bodies that moves periodically solving the three-body problem. The following table provides the initial conditions of these 9 periodic solutions. The first entry $T$ represent the period of the reduced periodic solution. This is, the positions and velocities after  $T$ units of time agree with the initial positions and velocities up to a rotation of $\theta$ radians. This rotation $\theta$ is given in the second column of the table. Every color in these images show the trajectory after $T$ units of time. These points are part of the path $S_1$. These periodic solutions are odd/even solutions and therefore they have two symmetries. 

\begin{center}
    \begin{tabular}{ | l | l | l |}
    \hline
    Coordinates:\hskip .4cm $(T,a,b)$ & $\Theta(T,a,b)$ & Image \\
    & &  location \\ \hline
    (13.451787406253272, 2.8889144618608733, 2.8514754095228265)&$\frac{11 \pi }{6}$ & (1,1)  \\ \hline
    (13.60916416956247, 2.8011806024112604, 2.8994928664781003)& $\frac{13 \pi }{7} $& (1,2)  \\ \hline
   (13.719567754906889, 2.737827736709047, 2.9324647582021814)&$ \frac{15 \pi }{8}$ & (1,3) \\ \hline 
     (13.801004797570164, 2.689923890770301, 2.9564685196206826)& $\frac{17 \pi }{9}$ & (2,1) \\
    \hline  (14.320649658996734, 2.344448198979306, 3.106748248260848) &$ 2\pi$ & (2,2)\\
    \hline (14.657003574778068, 2.0202122629047206, 3.212388918731313) & $\frac{17 \pi }{8}$& (2,3)\\
    \hline(14.686554119081652, 1.9783743950165875, 3.2235432541405697)&$ \frac{15 \pi }{7}$&(3,1) \\
    \hline (14.719531067694582, 1.9241180805387452, 3.2371600512994565) & $\frac{13 \pi }{6}$& (3,2) \\
    \hline 
    (14.754026030069339, 1.8509347720107878`, 3.254002265714927) & $\frac{11 \pi }{5}$& (3,3) \\
    \hline 
        \end{tabular}
\end{center}

\begin{figure}[hbtp]
\begin{center}\includegraphics[width=1.1\textwidth]{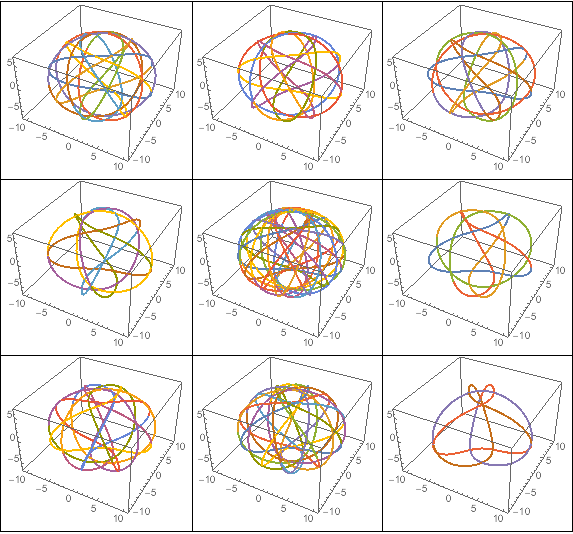}
\end{center}
\caption{Trajectory of one of the bodies for some of the periodic solutions described in this paper.}\label{t9}
\end{figure}

\subsection{Periodic solution displayed on Figure \ref{t9}}

Every entry of the matrix  in  Figure \ref{t9} (located at the end of the paper) shows the trajectory  of one of the three bodies that moves periodically solving the three-body problem. The following table provides the initial conditions of these 9 periodic solutions. The first entry $T$ represent the period of the reduced periodic solution. This is, the positions and velocities after  $T$ units of time agree with the initial positions and velocities up to a rotation of $\theta$ radians. This rotation $\theta$ is given in the second column of the table.  Every color in these images show the trajectory after $T$ units of time. These point are part of the path $S_1$. These periodic solutions are odd/even solutions and therefore they have two symmetries. 

\begin{center}
    \begin{tabular}{ | l | l | l |}
    \hline
    Coordinates:\hskip .4cm $(T,a,b)$ & $\Theta(T,a,b)$ & Image \\
    & &  location \\ \hline
    (14.782277611145048, 1.7467887078517095, 3.274928737204819)&$\frac{9 \pi }{4}$ & (1,1)  \\ \hline
    (14.786834114569135, 1.6762115129006632, 3.287061516889608)& $\frac{16 \pi }{7} $& (1,2)  \\ \hline
   (14.774959957278682, 1.5866395627104857, 3.300053407198923)&$ \frac{7 \pi }{3}$ & (1,3) \\ \hline 
     (14.7286222583901, 1.4691728012133716, 3.3129725743424996)& $\frac{12 \pi }{5}$ & (2,1) \\
    \hline  (14.699688133885457, 1.4214348721248662, 3.316877441047812) &$ \frac{17 \pi }{7}$ & (2,2)\\
    \hline (14.60791563192398, 1.3083245035642173, 3.3230110738911285) & $\frac{5 \pi }{2}$& (2,3)\\
    \hline(14.4942511626078, 1.2033604318255704, 3.324782465355579)&$ \frac{18 \pi }{7}$&(3,1) \\
    \hline (14.444161869664121, 1.163467447374621, 3.3244726259938155) & $\frac{13 \pi }{5}$& (3,2) \\
    \hline 
    (14.319901500300364, 1.0746785092812807, 3.321867812104126) & $\frac{8 \pi }{3}$& (3,3) \\
    \hline 
        \end{tabular}
\end{center}

\begin{figure}[hbtp]
\begin{center}\includegraphics[width=1.1\textwidth]{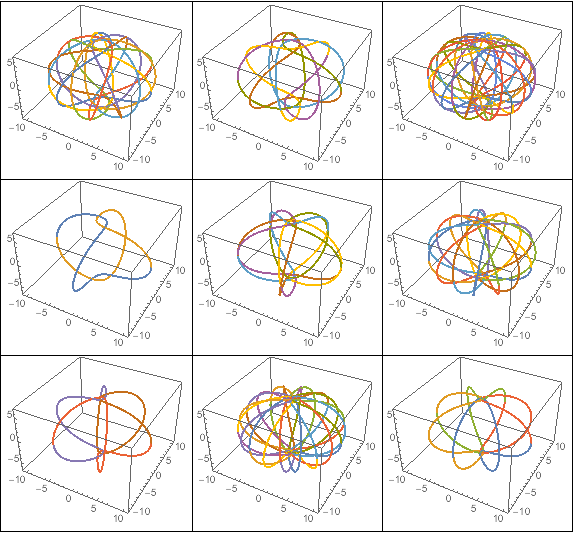}
\end{center}
\caption{Trajectory of one of the bodies for some of the periodic solutions described in this paper.}\label{fth9}
\end{figure}

\subsection{Periodic solution displayed on Figure \ref{fth9}}

Every entry of the matrix  in  Figure \ref{fth9} (located at the begging of the paper) shows the trajectory  of one of the three bodies that moves periodically solving the three-body problem. The following table provides the initial conditions of these 9 periodic solutions. The first entry $T$ represent the period of the reduced periodic solution. This is, the positions and velocities after  $T$ units of time agree with the initial positions and velocities up to a rotation of $\theta$ radians. This rotation $\theta$ is given in the second column of the table.  Every color in these images show the trajectory after $T$ units of time. These point are part of the path $S_1$. These periodic solutions are odd/even solutions and therefore they have two symmetries. The periodic solution given by the last row in this table is very close to a double collision.

\begin{center}
    \begin{tabular}{ | l | l | l |}
    \hline
    Coordinates:\hskip .4cm $(T,a,b)$ & $\Theta(T,a,b)$ & Image \\
    & &  location \\ \hline
    (14.155328130457113, 0.9714850334241202, 3.3156069432644326)&$\frac{11 \pi }{4}$ & (1,1)  \\ \hline
    (14.054387879818348, 0.9133674393862702, 3.310617589529648)& $\frac{14 \pi }{5} $& (1,2)  \\ \hline
   (13.986985590840762, 0.8760917450402825, 3.306890869738851)&$ \frac{17 \pi }{6}$ & (1,3) \\ \hline 
     (13.659851940872658, 0.70551977441241, 3.285155344523193)& $3 \pi $ & (2,1) \\
    \hline  (13.315322141213205, 0.5298843477532729, 3.2569218042195365) &$ \frac{16 \pi }{5}$ & (2,2)\\
    \hline (13.240214962690887, 0.4900089555434633, 3.250127637396059) & $\frac{13 \pi }{4}$& (2,3)\\
    \hline(13.125795380869723, 0.4265448292958072, 3.239357536186598)&$ \frac{10 \pi }{3}$&(3,1) \\
    \hline (13.0441592363555, 0.3781732129466524, 3.23136706052875) & $\frac{17 \pi }{5}$& (3,2) \\
    \hline 
    (12.938357527438528, 0.30903261362830825, 3.2206342443812517) & $\frac{7 \pi }{4}$& (3,3) \\
    \hline 
        \end{tabular}
\end{center}

 \vfil
 \eject
 
 \vfil
 \eject
 \vfil
 \eject
 \vfil
 \eject

\end{document}